\documentclass[12pt]{article}
\usepackage{graphicx}
\usepackage{amsmath}
\usepackage[utf8]{inputenc}
\usepackage{tikz}
\usepackage{tikz-3dplot}
\usepackage{url}
 \vspace{-8ex}
  \date{}

\begin{document}
\begin{titlepage}

\titlepage
\title{Map projection: A mathematical investigation on the distance-preserving property of an equidistant cylindrical  projection}
\author{Bingheng Yang}
\date{\url{yangbingheng00@gmail.com}}

\maketitle

\end{titlepage}

\newpage

\tableofcontents

\contentsline {section}{\numberline {1}Abstract}{2}
\contentsline {section}{\numberline {2}Introduction}{3}
\contentsline {subsection}{\numberline {2.1}A history of cartography}{3}
\contentsline {section}{\numberline {3}Our map}{4}
\contentsline {subsection}{\numberline {3.1}Equidistant cylindrical}{4}
\contentsline {subsection}{\numberline {3.2}Projection method}{6}
\contentsline {subsection}{\numberline {3.3}Forming the projection formulas}{8}
\contentsline {section}{\numberline {4}Distortion in distance}{11}
\contentsline {subsection}{\numberline {4.1}About distortion}{11}
\contentsline {subsection}{\numberline {4.2}How can the distortion be seen?}{12}
\contentsline {section}{\numberline {5}Investigation methods}{14}
\contentsline {subsection}{\numberline {5.1}Geometrical analysis}{14}
\contentsline {subsubsection}{\numberline {5.1.1} Tissot's indicatrix} {15}
\contentsline {subsubsection}{\numberline {5.1.2} Maximum angular and distance distortion } {17}
\contentsline {subsection}{\numberline {5.2}Trigonometrical analysis}{22}
\contentsline {subsubsection}{\numberline {5.2.1}Spherical coordinates}{23}
\contentsline {subsubsection}{\numberline {5.2.2}Relationship between straight-line distance and real distance}{25}
\contentsline {section}{\numberline {6}Conclusion}{28}
\contentsline {section}{\numberline {A}Appendix}{31}
\contentsline {subsection}{\numberline {A.1}Glossary}{31}
\contentsline {subsection}{\numberline {A.2} General model for equidistant cylindrical projection}{32}
\contentsline {subsection}{\numberline {A.2}The Haversine formula} {32}
\contentsline {subsection}{\numberline {A.3}Bibliography}{33}

\newpage

\section*{Abstract}

A normal world map can be defined as a representation on a flat surface that shows the features of the earth, such as the area of Australia or the latitude of the arctic circle, in their respective forms, sizes, and relationships according to some convention of representation. 

All world maps have inevitable errors\footnote{Vox, \textsl{Why all world maps are wrong? (2016, December 2)}

 https://www.youtube.com/watch?v=kIID5FDi2JQ}: Some maps fail to represent correct directions, some give us wrong displays of areas. In other words, there exists no such maps that could show all the right characteristics of the Earth as the globe would. This is due to the distortion caused by different projection methods.\footnote{Numberphile, \textsl{A strange map projection (Euler spiral) (2018, November 13)}

https://www.youtube.com/watch?v=D3tdW9l1690}

This research work aims to explore the distortions in distance in equidistant cylindrical projection. Although the projection is described as "equidistant", i.e. distance-preserving, it is far from error-free. By constructing a realistic and appropriate model of the projection, this work will demonstrate  \textbf{to what extend are distortions in distance present in equidistant cylindrical projection} using mathematical methods such as linear modeling, differentiation and trigonometric relationships. 

The investigation is conducted by examining the projection performing geometrical and trigonometrical analyses. The horizontal bending that occurs in the projection process can be assessed by performing a geometric analysis using Tissot's indicatrices. In addition, the concept of the spherical coordinates, alongside with trigonometrical identities, can be used to illustrate the route from a point to another as a curve. With a combination of the knowledge extracted from the examination of the projection using those two theories, this research aims to fully unravel the degree of distortion in distance in equidistant cylindrical projections.
\newpage

\begin{center}
\section*{Introduction}
\end{center}

\subsection*{A history of cartography}

Cartography\footnote{James S. Aber, \textsl{Brief history of maps and cartography}, retrieved March 2019.

\url{http://academic.emporia.edu/aberjame/map/h_map/h_map.htm}}, or the study of making maps, is one of the oldest research concepts that human practiced. For thousands of years, mathematicians and astronomers have been successful to accumulate world maps that accurately display properties of Earth to scale. The idea of maps already existed during the Babylonian era around 2300 B.C, and the Greek developed world map designs to an advanced level. 

Later, during the end of the Middle Ages, the first concept of a world map, the so called T-O map, was invented. For the next few hundred years, exploration trips, improvement in technology and the increasing market need lead to a rapid development phase of world maps. Different designs were created by talented cartographers for a number of purposes. Some of the world map projections are for example, the Mercator projection and the Eckert projections. 

After Soviet engineers launched the first  satellite Sputnik 1, the design process of world map projections became much easier. The rise of GIS - Geographic Information Systems - enabled the creation of more accurate maps and newer concepts such as digitalized maps, e.g. Google maps.

\newpage
\section*{Our map}

In our projection, we first assume the Earth to be \textbf{spherical} and the radius of the Earth to be \textbf{6371 km}. We aim to represent a projection using Cartesian xy-coordinates, with the x-axis describing the East-West direction and the y-axis demonstrating the North-South direction. The scale factors, which indicate how does the real distance between two points on the surface of the Earth correspond to the distance displayed by the projection, are also represented by mathematical formulas.

We must be aware of the properties and features of the projections. In the following sections, the formation method of the map projection is investigated, from which the formulas of an equidistant cylindrical projection are found. All numerical values for angles will be in \textbf{radians}. All numerical values for distances will be in \textbf{kilometers}.

\subsection*{Equidistant cylindrical}

Equidistant cylindrical projection\footnote{Wikipedia, \textsl{Equidistant cylindrical projection,} retrieved December 2018.} is one of the first world map projections created. It was first established by Claudius Ptolemy, and improved by later generations. The purpose of the projection is to give an accurate presentation of the distance between two points on the surface of the Earth: By multiplying the distance between two points on the map with a scale factor, the result should correspond to the real distance.

A \textbf{cylindrical projection} is created by first \textbf{inscribing} a sphere into a right cylinder without the top and bottom surfaces, then transferring the points on the spherical grid to the surface of the cylinder. Finally, the side surface of the cylinder is unfolded, resulting a rectangular projection. We can see the construction process\footnote{Eotvos Lorand University, \textsl{Cylindrical projections}, retrieved July 2019.

\url{http://lazarus.elte.hu/cet/modules/guszlev/cylin.htm}} in figure 2. This will give us a rectangular-shaped map, such as in figure 1: In equidistant cylindrical projection, the shape of the world is rectangular. By name, the projection should be distance-preserving, meaning that the distance on the map corresponds to the real distance on the globe. From figure 2, we can observe a number of properties that are unique for equidistant cylindrical projection: All vertical meridians are parallel to each other and all horizontal parallels are also parallel to each other. The gap between individual meridians and the gap between individual parallels are all equal in the projection. Standard parallels and meridians are all straight lines which are perpendicular to each other: The square grids, which are formed by the intersections of meridians and parallels are also all equal in size. The length of meridians is half the length of parallels.

\begin{figure}
\includegraphics[scale=0.9]{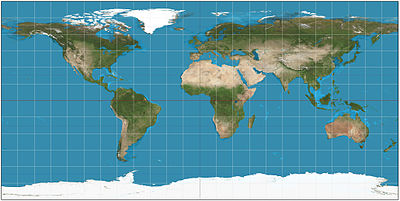}
\caption{An example of an equidistant cylindrical projection.}
\end{figure}

\begin{figure}

\begin{tikzpicture}

\draw (0,-1,-3)  ellipse (1 and 0.5);
\draw (1.1,-2,0) circle [radius=1];
\draw (-1,-1,-3)--(0,-4,0);
\draw (1,-1,-3)--(2.1,-4,0);
\draw (1.05,-4,0) ellipse (1.05 and 0.5);

\draw [dashed] (0.1,-2,0) to [in=25,out=55] (2,-2,0);
\draw [dashed] (0.1,-2,0) to [in=-115,out=-115] (2,-2,0);
\draw [dashed] (1,-3,0) to [in=0,out=0] (1,-1,0);

\draw [->] (4,-2,0) -- (6,-2,0)[anchor=right];
\draw (8,-1,0)--(12,-1,0)--(12,-3,0)--(8,-3,0)--(8,-1,0);
\draw [dashed] (9,-1,0)--(9,-3,0) node[below]{M};
\draw [dashed] (8,-1.5,0)--(12,-1.5,0) node[right]{P};
\draw [dashed] (10,-1,0)--(10,-3,0);
\draw [dashed] (11,-1,0)--(11,-3,0);
\draw [dashed] (8,-2,0)--(12,-2,0);
\draw [dashed] (8,-2.5,0)--(12,-2.5,0);
\draw [->] (-0.35,-2,0)--(1.5,-2,0);
\node at (-0.5,-2,0) {$M$};
\draw [->] (3,-1.5,0)--(1,-1.6,0);
\node at (3.2,-1.5,0) {$P$};

\end{tikzpicture}
\caption{The sphere is placed in a topless and bottomless right circular cylinder. Individual points on the sphere are transferred onto the side of the cylinder; An image of a rectangle is generated after the cylinder is unfolded, because the side surface of a right cylinder demonstrates a rectangle. Meridians (marked with M) and standard parallels (marked with P) intersect each other at right angles and they are equally spaced. $|M|=\frac{1}{2}|P|$ holds for all meridians and parallels in both the globe and the projection.}
\label{Sphere in cylinder}
\end{figure}
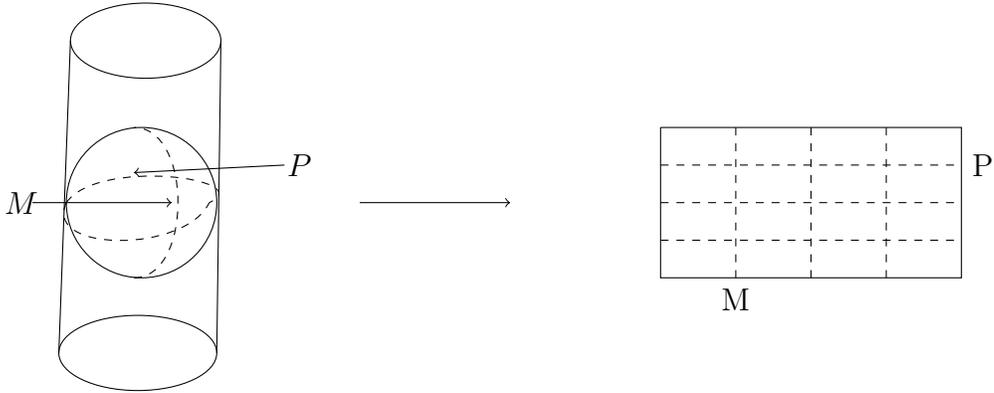

\newpage

\subsection*{Projection method}

In this section, the way how individual points are transferred onto the cylinder is explained and justified. 

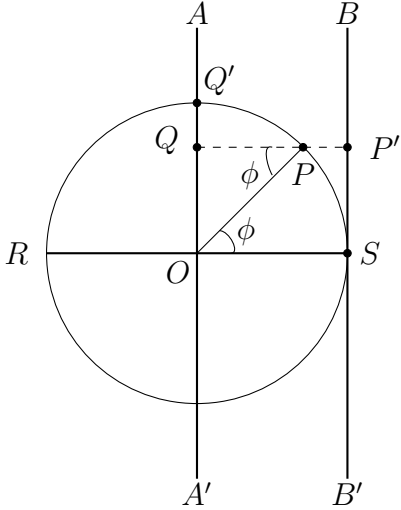
\begin{figure}
\begin{tikzpicture}
\draw (0,0) circle [radius=2];
\draw [thick] (0,-3) -- (0,3);
\draw [thick] (-2,0) -- (2,0);
\draw (0,0) -- (1.41,1.41);
\draw [thick] (2,-3) -- (2,3);
\draw (0.45,0) to [out=0,in=0] (0.3,0.3);
\draw [dashed] (0,1.41) -- (1.41,1.41);
\draw (1,1.41) to [out=160,in=120] (1,1.04);
\node at (-0.25,-0.25) {$O$};
\node at (0,3.2) {$A$};
\node at (2,3.2) {$B$};
\node at (2.3,0) {$S$};
\node at (0.65,0.3) {$\phi$};
\node at (-0.4,1.5) {$Q$};
\node at (1.41,1.05) {$P$};
\node at (2.5,1.41) {$P'$};
\node at (-2.4,0) {$R$};
\node at (0,-3.2) {$A'$};
\node at (2,-3.2) {$B'$};
\node at (0.7,1.04) {$\phi$};
\node at (0.3,2.3) {$Q'$};
\draw [fill] (2,0) circle [radius=0.05];
\draw [fill] (0,1.41) circle [radius=0.05];
\draw [fill] (2,1.41) circle [radius=0.05];
\draw [fill] (1.41,1.41) circle [radius=0.05];
\draw [fill] (0,2) circle [radius=0.05];
\draw [dashed] (1.41,1.41) -- (2,1.41);
\end{tikzpicture}
\caption{A front view of figure 2. The angle $\phi$ is preserved in both 2D and 3D illustrations.}
\label{Front view of 3D figure}
\end{figure}

Figure \ref{Sphere in cylinder} implies that the sides of the cylinder are \textbf{tangents} to the sphere. To demonstrate this observation, figure \ref{Front view of 3D figure} shows the front view of figure \ref{Sphere in cylinder}. The circle with diameter $RS$ - which indicates the sphere - has a tangent $BB'$, which denotes the side of the cylinder.

In figure \ref{Front view of 3D figure}, $O$ denotes the middle point of the sphere, and let  points $S$ and $P$ be on the sphere. $AA'$ is a line crossing the diameter of the circle. It is said in earlier sections that points $S$ and $P$ have to be transformed onto the side of the cylinder, i.e. tangent $BB'$. Since the circle and $BB'$ touch at point $S$, $S$ is on the side of the cylinder. However, in the case of point $P$, it is projected to $P'$.

Another point $Q$ on the sphere (circle) is selected such that $Q$ has the same latitude as $P$ and $Q$ is always located on $AA'$. Lines $OP$ and $OS$ form an angle $\phi$. When point $P$ travels anti-clockwise to approach line $AA'$, $Q$ translate upwards and meets $P$ at $Q'$, during when the position of $P'$ on $BB'$ also shifts upwards. Additionally, as $|PP'|$ increases, angle $\phi$ increases when $P$ approaches $Q'$.

\begin{figure}
\includegraphics[scale=0.5]{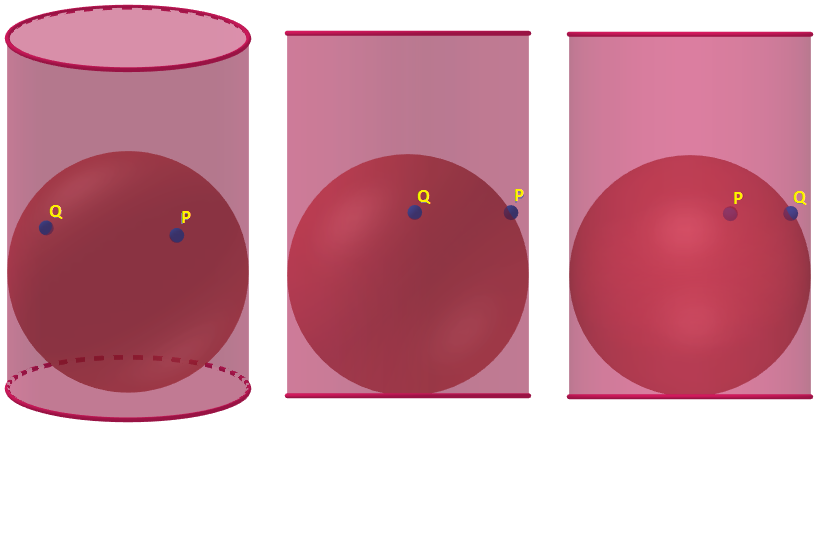}
\caption{The sphere and cylinder in figure \ref{Sphere in cylinder} from different perspectives. As the third cylinder illustrates, the sphere can be rotated so that point $Q$ can be on a position equivalent to the position of point $P$ in the second cylinder. Hence, the position of point $P$ in figure \ref{Front view of 3D figure} can be recreated for all points along the same longitude, such as point $Q$. Hence, the points $P$ and $Q$ undergo the same extent of disposition, causing the Euclidean distance $|PQ|$ (figure \ref{Front view of 3D figure}) to distort. Since $P$ and $Q$ share the same longitude, it seems that distortion in distance is independent from the longitudinal coordinates.}
\label{rotation}
\end{figure}

From figure \ref{Front view of 3D figure}, the scale factors can be deduced. As figure \ref{rotation} shows, the sphere can be rotated in a way so that point $Q$ is transformed to the position of point $P$ in figure \ref{Front view of 3D figure}, and point $Q$ will be projected similarly to the projection of point 
$P$ (page \pageref{Front view of 3D figure}). Thus,  distortion in distance between two points on the projection is independent of their horizontal coordinates. Therefore, only vertical translation can affect the distance-preserving ability of the projection. 
 
Figure \ref{Sphere in cylinder} also shows that a rectangular map projection is formed. Let the horizontal side of the projection be the length and the vertical side of the projection be the width. Since the side surface of the cylinder is a tangent to the sphere, the length of the projection should correspond to the perimeter of the circle in figure \ref{Front view of 3D figure}: if we consider $|RS|$ to be on the sphere (figure \ref{Front view of 3D figure}), it corresponds to an half-circle that has the radius of the original sphere; the half-circle would always be on the cylinder, hence the length of the cylinder equals the perimeter of the half-circle in figure \ref{Front view of 3D figure}. Similarly, the width of the projection should be set to include all projected points such as $P'$: when the radius of the circle (figure \ref{Front view of 3D figure}) increases, the position of $P'$ moves up in relative terms. Hence, both the length and the width of the projection are dependent on the radius of the sphere.

\subsection*{Forming the projection formulas}

Using the information acquired in the last section, we strive to formulate appropriate formulas for an equidistant cylindrical projection using the Cartesian grid. The points of longitude are denoted by $\lambda$ and the points of latitude are denoted by $\phi$. Let the relative scale factor of the meridian be $h$ and the relative scale factor of the standard parallel be $k$.

\begin{figure}
\begin{tikzpicture}

\draw (0,0) -- (-1.5,0);
\draw (0,0) -- (0,3);
\draw [dashed] (0,0) -- (0,-3);
\draw [dashed] (0,0) -- (1.8,-2.4);
\draw (0,0) circle [radius=3];
\draw (0,0) circle [radius=1.5];
\draw [fill] (0,-3) circle [radius=0.05];
\draw [fill] (0,0) circle [radius=0.05];
\draw [fill] (0,-1.5) circle [radius=0.05];
\draw [fill] (1.8,-2.4) circle [radius=0.05];
\draw [fill] (0.9,-1.2) circle [radius=0.05];
\draw [dashed] (0,-0.5) to [out=0,in=0] (0.2,-0.3);
\node at (-0.75,0.25) {$r$};
\node at (0.25, 2.25) {$R$};
\node at (-0.3,-1.8) {$Q$};
\node at (-0.3,-2.6) {$Q'$};
\node at (0.5,-1.2) {$P$};
\node at (1.9,-1.9) {$P'$};
\node at (0.4,0) {$\phi$};
\end{tikzpicture}
\caption{Above view of figure \ref{rotation}. The greater circle stands for the circular base of the right cylinder, and the smaller circle illustrates the above image of the inscribed sphere. Points $P$ and $Q$ are located on a circle with radius $r$, and the circular base of the cylinder has a radius $R$. $\phi$ is the common latitude of all points $P$ and $Q$ on the sphere and $P'$ and $Q'$ on the projection.}
\label{above}
\end{figure}
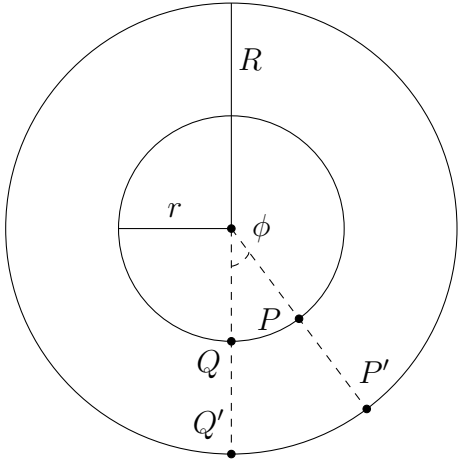

As mentioned earlier, horizontal translations will not cause distortion in distance. Hence, $h$ is a constant. To make further calculations simpler, we can let $h=1$. On the other hand, to calculate $k$, the relationship between $|PQ|$ and $|P'Q'|$ must be known; this can be deduced from figure \ref{above}. Because the scale factor\footnote{See glossary for definition.} is calculated as $$\frac{\text{Euclidean distance between two points on the projection}}{\text{Euclidean distance between two points on the surface of the Earth}},$$ we deduce that $k=\dfrac{|P'Q'|}{|PQ|}$ from figure \ref{above}. Also, $|PQ|=r\phi$ and $|P'Q'|=R\phi$. Thus, 

$$k=\frac{|P'Q'|}{|PQ|}=\frac{R\phi}{r\phi}=\frac{R}{r}.$$

Notice that the relationship between $R$ and $r$ can be determined using figure \ref{Front view of 3D figure}. Since $R$ is the radius of the base circle of the cylinder, it corresponds to $|OS|$, which equals $|OP|$. Similarly, $r$ stands for $PQ$. Using trigonometrical identities, we deduce $r=R\cos\phi$. Therefore, 

$$k=\frac{R}{r}=\frac{R}{R\cos\phi}=\frac{1}{\cos\phi}.$$

Hence, the scale factors of the projection are 

\begin{equation}
h=1
\label{eq.1}
\end{equation}

\begin{equation}
k=\frac{1}{\cos\phi}
\label{eq.2}
\end{equation}

Next, a rectangular projection is built. Since the size of the projection depends on the radius of the sphere (figure \ref{Sphere in cylinder}), let the radius of the sphere be $R$. 

One way to create the model is to define the intersection of the central meridian and standard parallel to be at the origin. By this action, \textsl{the central meridian is chosen to be the Greenwich line ($\lambda=0$) and the standard parallel is chosen to be the equator ($\phi=0$)}. Consequently, the x-coordinate is dependent on $\lambda$ and the y-coordinate on $\phi$. Recall from figure \ref{Sphere in cylinder} that all meridians are equal in length and all parallels are equal in length; Furthermore, the meridians are a half of the sphere's outer perimeter and the parallels are equal to the sphere's outer perimeter. Hence, we deduce that the length of each individual meridian must be $\pi{R}$ and the length of each individual parallel must be $2\pi{R}$, so the length of the projection is going to be $2R\pi$ and the width of the projection is going to be $R\pi$. The latitude coordinates of the Earth varies between $90^{\circ}$N and $90^{\circ}$S and the longitude coordinates shifts from $180^{\circ}$E to $180^{\circ}$W; To make modeling easier, we let latitude points south of the standard parallel and longitude points west of the prime meridian to take \textbf{negative values} in radians. Therefore, the coordinates of a random point  located on the projection satisfy 

\begin{equation}
y=R\phi, -\frac{\pi}{2}\leq\phi\leq\frac{\pi}{2}
\label{eq.3}
\end{equation}

\begin{equation}
x=R\lambda, -\pi\leq\lambda\leq\pi.
\label{eq.4}
\end{equation}

In this way, $\phi_{\text{max}}=\pi$ and $\lambda_{\text{max}}=2\pi$, so the length of the meridians and parallels will be $\pi{R}$ and $2\pi{R}$. Hence, by combining equations \ref{eq.1}, \ref{eq.2}, \ref{eq.3} and \ref{eq.4}, the equations for equidistant cylindrical projection are:

$$y=R\phi$$

$$x=R\lambda$$

$$h=1$$

$$k=\dfrac{1}{\text{cos}\phi}$$

$$-\frac{\pi}{2}\leq\phi\leq\frac{\pi}{2}, -\pi\leq\lambda\leq\pi.$$

Where 
\begin{itemize}

\item{$x$ is the coordinates of East-West axis of the map}

\item{$y$ is the coordinates of North-South axis of the map}

\item{$h$ is the relative scale factor of the meridian}

\item{$k$ is the relative scale factor of the standard parallel}

\item{$R$ is the radius of the Earth}

\item{$\phi$ is latitude in radians}

\item{$\lambda$ is longitude in radians.}

\end{itemize}
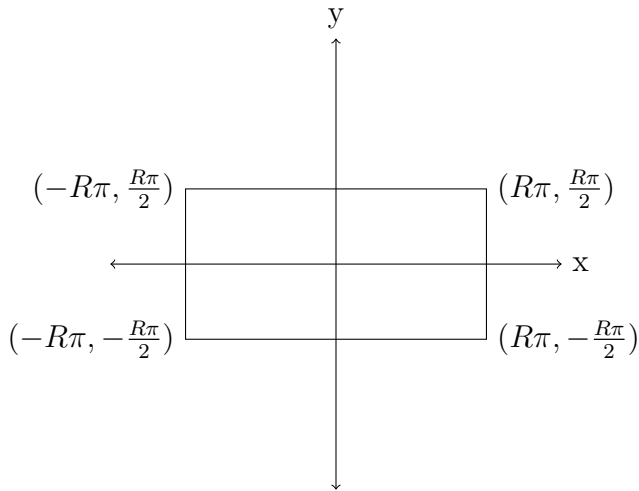
\begin{figure}

\begin{tikzpicture}[x=0.1cm,y=0.1cm]

\draw [<->] (-30,0)--(30,0) node[right]{x} ;
\draw [<->] (0,-30)--(0,30) node[above]{y} ;
\draw (-20,-10) rectangle (20,10);

\node [left] at (-20,-10) {($-R\pi,-\frac{R\pi}{2}$)};
\node [right] at (20,-10) {($R\pi,-\frac{R\pi}{2}$)};
\node [left] at (-20,10) {($-R\pi, \frac{R\pi}{2})$};
\node [right] at (20,10) {($R\pi, \frac{R\pi}{2}$)};

\end{tikzpicture}
\caption{Equidistant cylindrical projection in xy-coordinates.}
\label{Formulas}
\end{figure}
As noted earlier, in this particular model, the prime meridian is set to be the Greenwich line, which has the longitude of $0^\circ$, corresponding to the y-axis. Contrarily, the standard parallel is selected to be the equator. All points along the Equator have $0^\circ$ of latitude, therefore aligning to the x-axis. A model of the projection is shown by figure \ref{Formulas}.

From the formulas we can already see one possible cause of distortion: $k$, the relative scale factor of the horizontal parallels, is a variable. This means that the distance between two points with different longitude will not always correspond to the real distance. When $k=1$, $\cos\phi=0 \leftrightarrow \phi=0$, meaning that the error in distance is minimized at the equator. The distance distortion is at the maximum in the poles, since $\displaystyle{\lim_{\phi\to\frac{\pi}{2}-} \frac{1}{\cos\phi}} = \infty$ and $\displaystyle{\lim_{\phi\to-\frac{\pi}{2}+} \frac{1}{\cos\phi}} = \infty$.
\section*{Distortion in distance}

\subsection*{About distortion}

\footnote{Rice University, \textsl{Mapping the sphere}, retrieved April 2019.

\url{https://math.rice.edu/~polking/cartography/cart.pdf}} Distortion is and will always be present in map projections. There exists no maps that can perfectly preserve and represent the most important geometrical identities of the Earth to scale: the latitude and longitude angles, the surface area of regions, and the distance between two points on the earth. The distortion will be or will not be significant depending on the purpose of the map. For example, if we want to use a map for sailing from Gibraltar to Cape Town, we can ignore the irrelevant area distortion as long as the map is conformal so the direction angle is always correct. On the other hand, the distortion in distance in a "distance-preserving" projection must be taken into account. In the following section, the distance distortion in equidistant cylindrical projection is demonstrated.

\subsection*{How can the distortion be seen?}
\begin{figure}
\begin{tikzpicture}[x=0.1cm,y=0.1cm]

\draw [<->] (-30,0)--(30,0) node[right]{$x$} ;
\draw [<->] (0,-30)--(0,30) node[above]{$y$} ;
\draw (-20,-10) rectangle (20,10);
\draw [fill] (4,4.05) circle [radius=0.5] node[right] {P$_1$};
\draw [fill] (-4.2,-1.85) circle [radius=0.5] node[left] {P$_2$};
\draw (4,4.05)--(-4.2,-1.85);
\node at (-3,2) {$d$};
\end{tikzpicture}
\caption{P$_1$ and P$_2$ on a modeled equidistant cylindrical projection. 
$d$ demonstrates the distance between P$_1$ and P$_2$.}
\label{P1 and P2}
\end{figure}
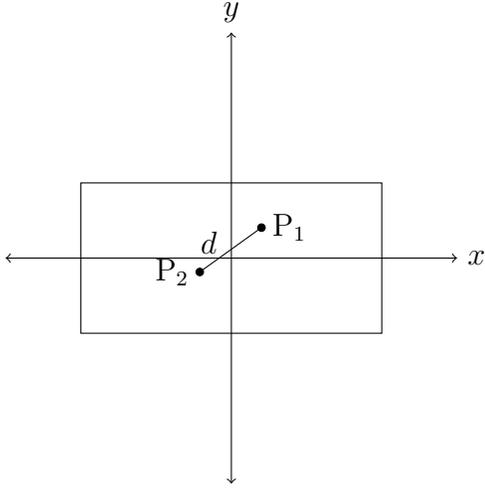

Suppose we want to determine the distance between two points, P$_1$ and P$_2$, using the formulas of equidistant cylindrical projection. According to the formulas (page 7), by recognizing the longitude and latitude coordinates of any two points, we can deduce their x- and y-coordinates. Therefore, we can arbitrarily select the coordinates of P$_1$ and P$_2$: For example, let P$_1=$ ($24.3^\circ$E, $23.4^\circ$N) and P$_2=$ ($39.2^\circ$W, $3.67^\circ$S). Their relative distance on the map is illustrated in figure \ref{P1 and P2}. To calculate the distance between P$_1$ and P$_2$, we need to first change the longitude and latitude in the coordinates into radians.

P$_1=(24.3^\circ$E, $23.4^\circ$N) $ = (0.424...,0.408...)$

P$_2=(39.2^\circ$W, $3.67^\circ$S) $ = (-0.684...,-0.064...)$

By using the formulas of equidistant cylindrical projection (page 7), we can calculate the values of x and y-coordinates.

$x_1=R\lambda_1= 6371\times0.424...\approx2702.036665$

$x_2=R\lambda_2= 6371\times(-0.684...)\approx-4358.841336$

$y_1=R\phi_1=6371\times0.408...\approx2601.960997$

$y_2=R\phi_2=6371\times(-0.064)\approx-408.0853582$

The distance between any two points in the xy-plane can be determined by the relationship $d=\sqrt{{(x_1-x_2)}^2+{(y_1-y_2)}^2}$. Hence, the distance between P$_1$ and P$_2$ is 

$$d=\sqrt{(2702...-(-4358...))^2+(2601...-(-408...))^2}=7675.700438\approx7676.$$

However, by using \textsl{Movable type script}\footnote{Movable type scripts, \textsl{Calculate distance, bearing and more between latitude/longitude points}

https://www.movable-type.co.uk/scripts/latlong/html}, an engine that is able to generate real distances between two points on the surface of Earth, we observe that the real distance between our randomly selected P$_1$ and P$_2$ is 7502 km. Thus, our numerically calculated result has an error of $\dfrac{7676-7502}{7502}\times 100\%\approx2.32\%$. Therefore, we conclude that some distortion in distance appears in equidistant cylindrical projection. One possible reason of distortion could be that the Earth is assumed to be a perfect sphere, although the Earth is in fact, an ellipse in the real world. Nonetheless, in the following sections of this research, other possible sources of error are explored and evaluated. 

\newpage
\section*{Investigation methods}

In the previous section, we have demonstrated the existence of distance distortion in a equidistant cylindrical projection. We will investigate this issue by performing a \emph{geometrical} analysis and a \emph{trigonometrical analysis} to scrutinize the distance distortions of equidistant cylindrical projection. In this way, we aim to construct general methods to determine the extent of error in distance.

\subsection*{Geometrical analysis}

In previous sections, we have shown that $k$, the relative scale factor of the horizontal parallels, is a variable; the change in distortion is dependent on the value of latitude. Recall from page 9 that $\displaystyle{\lim_{\phi\to\frac{\pi}{2}-} \frac{1}{\cos\phi}} = \infty$ and $\displaystyle{\lim_{\phi\to-\frac{\pi}{2}+} \frac{1}{\cos\phi}} = \infty$: the horizontal distortion will increase out of bounds in polar areas.
\begin{figure}
\includegraphics[scale=0.9]{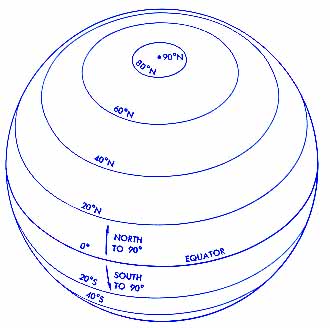}
\caption{Parallels on a map. It can be clearly seen that they are not equal in length.}
\label{Standard Parallels = circles}
\end{figure}

But why is this the case? Going back to figure 2, we notice that the projection is constructed such that all parallels are assumed to have the same length; The mapmaker chooses a standard parallel, and other parallels will be presumed to have the same length. However, from figure \ref{Standard Parallels = circles}, it can be observed that it is not the case. The length of parallels (except the standard parallel) are artificially changed, and distortion in distance is already inevitable. Therefore, we can deduce that artificial bending, also called as \textbf{flexion}\footnote{David M. Goldberg, J. Richard Gott III (2007), \textsl{Flexion and Skewness in Map Projections of the Earth}

http://www.physics.drexel.edu/~goldberg/projections/goldberggott.pdf}, occurs in equidistant cylindrical projection, which is due to the nature of its construction process. 

\begin{figure}
\includegraphics[scale=0.5]{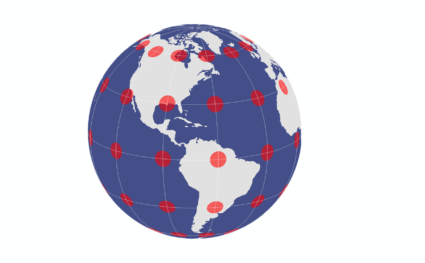}
\caption{Tissot's indicatrices in a sphere.}
\label{sphere}
\end{figure}

\begin{figure}
\includegraphics[scale=0.2]{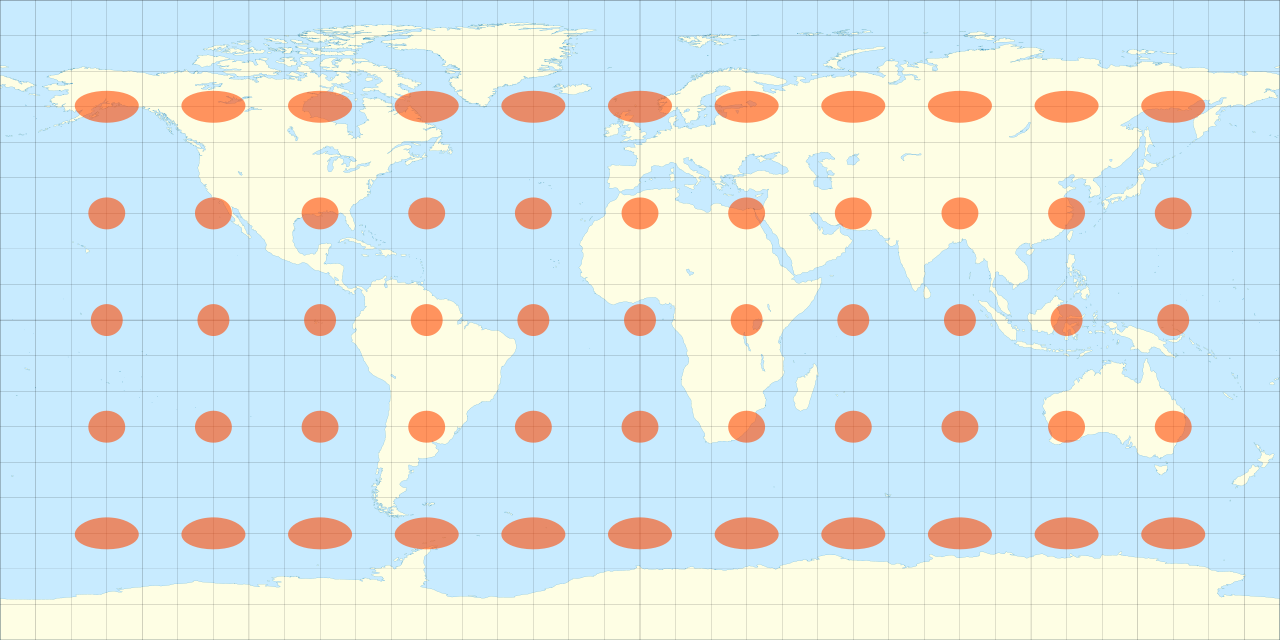}
\caption{Tissot's indicatrices in equidistance cylindrical projection. Notice the distortion at polar areas illustrated by the bending of indicatrices; the indicatrices become ellipses due to horizontal flexion.}
\label{indicatrix}
\end{figure}

\subsubsection*{Tissot's indicatrix}

There exists a various of methods to investigate the flexion which occurs in equidistant cylindrical projections. In this research, \textbf{Tissot's indicatrices}\footnote{John P. Snyder (1987), \textsl{Map projections - A working manual}} are used to illustrate the distortion in different regions. Tissot's indicatrices are circles with a tiny radius. They are constructed at given locations on a sphere, (see figure \ref{sphere}) and then the sphere undergoes cylindrical projection. As figure \ref{indicatrix} illustrates, the spheres closer to polar areas experience the most horizontal flexion and become elliptical. This observation corresponds to our formula for horizontal scale factor (page 7). Next, we examine the indicatrix in more detail and figure out to what extend does flexion cause distortions to appear.
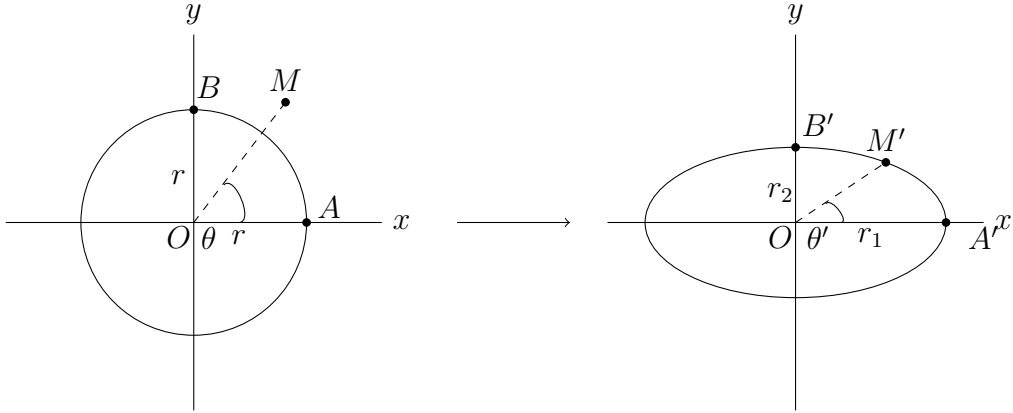
\begin{figure}

\begin{tikzpicture}[x=0.1cm,y=0.1cm]

\draw (-25,0)--(25,0) node[right] {$x$};
\draw (0,-25)--(0,25) node[above] {$y$};
\draw (0,0) circle [radius=15];

\draw [fill] (15,0) circle [radius=0.5];
\draw [fill] (0,15) circle [radius=0.5];
\draw [fill] (12.2,16) circle [radius=0.5] node[above] {$M$};
\draw [dashed] (0,0)--(12.2,16);
\draw (6,0) to [out=0,in=45] (4,5);
\node at (-2,-2) {$O$};
\node at (18,2) {$A$};
\node at (2,18) {$B$};
\node at (-2,6) {$r$};
\node at (6,-2) {$r$};
\node at (2,-2) {$\theta$};

\draw [->] (35,0)--(50,0);

\draw (55,0)--(105,0) node[right] {$x$};
\draw (80,-25)--(80,25) node[above] {$y$};
\draw (80,0) ellipse (20 and 10);
\draw [fill] (100,0) circle [radius=0.5];
\draw [fill] (80,10) circle [radius=0.5];
\draw [fill] (92,8) circle [radius=0.5] node[above] {$M'$};
\draw [dashed] (80,0)--(92,8);
\draw (86,0) to [out=0,in=45] (84,2.5);
\node at (105,-2) {$A'$};
\node at (83,13) {$B'$};
\node at (78,-2) {$O$};
\node at (78,4) {$r_2$};
\node at (90,-2) {$r_1$};
\node at (83,-2) {$\theta'$};
\end{tikzpicture}
\caption{An illustration of the extend of distortion of indicatrices caused by flexion due to the faulty nature of equidistant cylindrical projection. A small circle with radius $r$ experiences flexion and becomes an ellipse with a width $r_1$ and a height $r_2$.}
\label{Mashup}
\end{figure}

To clarify the transformation of the spheres in figure 8 to ellipses in figure 9, the process is modeled, as shown in figure 10. In \footnote{Wenping Jiang, Jin Li (2014), \textsl{The Effects of Spatial Reference Systems on the Predictive Accuracy of Spatial Interpolation Methods}

\url{https://www.researchgate.net/publication/259848543_The_Effects_of_Spatial_Reference_Systems_on_the_Predictive_Accuracy_of_Spatial_Interpolation_Methods}}
figure 10, we consider a small unit circle to be drawn on an arbitrary part of the surface of the Earth. By considering the circle to be small, we wish to generate an overall picture for distortions in distance between any two points on the projection, such as in figure \ref{indicatrix}: Greater bending shows more significant distortion in distance. The sphere is projected according to the methods of an equidistant cylindrical projection (see page 5 and figure 2). As shown in figure \ref{Front view of 3D figure}, in addition with the relationship $k=\frac{1}{\cos\phi}$, it can be deduced that the indicatrix is bent in the horizontal direction, which can be seen as a stretch from $OA$ to $OA'$ in figure \ref{Mashup}. Consequently, a vertical stretch also occurs, resulting $OB$ to decreases to $OB'$. As a result, the original circle is deformed and becomes an ellipse ($r\neq{r_1}\neq{r_2}$). 

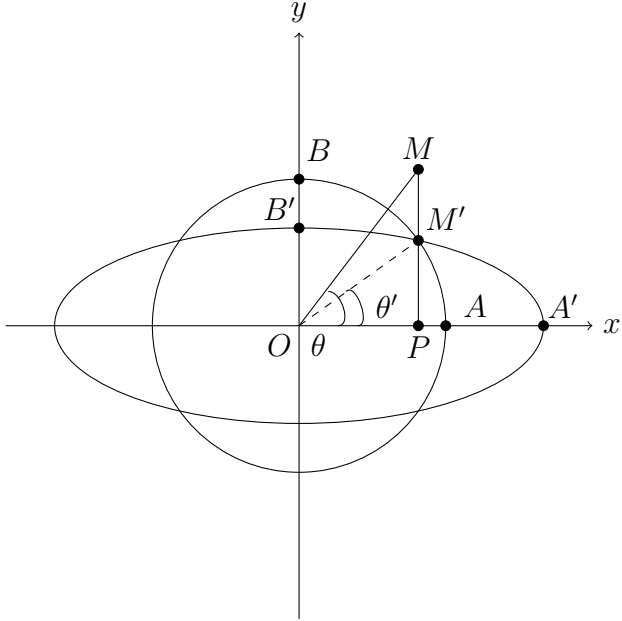
\begin{figure}

\begin{tikzpicture}[x=0.13cm,y=0.13cm]

\draw [->] (-30,0)--(30,0) node[right]{$x$} ;
\draw [->] (0,-30)--(0,30) node[above]{$y$} ;
\draw (0,0) circle [radius=15];
\draw (0,0) ellipse (25 and 10);
\draw (0,0)--(12.2, 16)--(12.2,0);
\draw [fill] (15,0) circle [radius=0.5];
\draw [fill] (0,15) circle [radius=0.5];
\draw [fill] (12.2,16) circle [radius=0.5] node[above] {$M$};
\draw [fill] (25,0) circle [radius=0.5];
\draw [fill] (0,10) circle [radius=0.5];
\draw [fill] (12.2,8.73) circle [radius=0.5];
\draw [fill] (12.2,0) circle [radius=0.5] node[below] {$P$};
\draw [dashed] (0,0)--(12.2,8.73);
\draw (4,0) to [out=0,in=0] (3,3.5);
\draw (6,0) to [out=0,in=45] (4.8,3.5); 
\node at (-2,-2) {$O$};
\node at (18,2) {$A$};
\node at (2,18) {$B$};
\node at (27,2) {$A'$};
\node at (-2,12) {$B'$};
\node at (15,11) {$M'$};
\node at (2,-2) {$\theta$};
\node at (9,2) {$\theta'$};

\end{tikzpicture}
\caption{A combined version of the bending of an indicatrix in figure \ref{Mashup}.}
\label{Extramashup}
\end{figure}

Consider point $M$ to be on the Earth's surface. Due to horizontal flexion, $M$ translates to $M'$ that is the point of intersection of lines $[MP]$ and the ellipse. In addition, angle $\theta$ decreases as $\Delta{OMA}$ is deformed and becomes $\Delta{OM'A}$. This phenomenon, know as an angular distortion, will affect the distance between two points on the projection. As a result, $|OM|\neq|OM'|$, therefore the distance between points $O$ and $M$ is not preserved. Therefore, the distortion in distance is $||OM|-|OM'||$. The overall combination of the two diagrams of figure \ref{Mashup} can be summarized using figure \ref{Extramashup}.

\newpage
\subsubsection*{Maximum angular and distance distortion}

In this section, the maximum angular distortion will be calculated and evaluated. Moreover, the horizontal distortion in distance due to a not-up-to-scale scale factor is analyzed.

Consider $\omega=\theta-\theta'$ where $\omega$ is the angular distortion. As can be seen in figure 11, the distortion is dependent on the coordinates of the point $M'$, which is projected according to the scale factors of the projection (see page 11).  Recall from page 8 that

$$y=R\phi$$

$$x=R\lambda$$

$$h=1$$

$$k=\frac{1}{\cos\phi}$$

\begin{figure}
\begin{tikzpicture}[x=0.13cm,y=0.13cm]

\draw (-10,0)--(20,0);
\draw (20,0)--(20,22);
\draw (20,22)--(-10,0);
\draw (-10,0)--(20,13);
\draw (20,0)--(18,0)--(18,2)--(20,2)--(20,0);
\draw [fill] (-10,0) circle [radius=0.5] node[left] {$O$};
\draw [fill] (20,13) circle [radius=0.5] node[right] {$M'$};
\draw [fill] (20,22) circle [radius=0.5] node[above] {$M$};
\draw [fill] (20,0) circle [radius=0.5] node[right] {$P$};
\draw (-6.2,0) to [out=0,in=0] (-6.2,3);
\draw (-4.2,0) to [out=0,in=0] (-4.2,2.4);
\draw (-3.2,2.9) to [out=0,in=0] (-3.2,5);
\node at (-8,-2) {$\theta$};
\node at (0,1.7) {$\theta'$};
\node at (-0.1,5.6) {$\omega$};
\node at (5,-4) {$\frac{1}{\cos\phi}$};
\node at (22,6) {$a$};
\node at (24,18) {$1-a$};
\end{tikzpicture}
\caption{The relationship between $\theta$, $\theta'$ and $\omega$ in relation to the changes of the value $a$.} 
\label{maximumomega}
\end{figure}
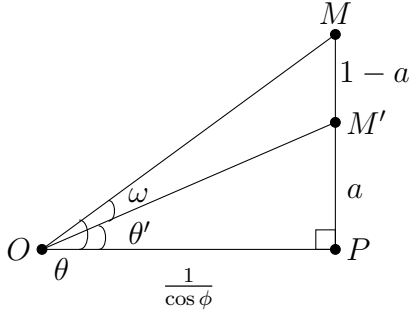

Thus, the latitude $\phi$ must be taken into account when computing the maximum angular distortion.  Figure \ref{maximumomega} represents the triangle $OMP$ in figure 11. In figure \ref{maximumomega}, it can be seen that $\vec{OM}=\vec{OP}+\vec{PM}$; the point $OM$ is transformed horizontally by 1 and vertically by 1. However, those transformations do not necessarily correspond to a distance of 1 unit on the map projection; the scale factors must be taken into account. Since the scale factors are $h=1$ and $k=\frac{1}{\cos\phi}$ respectively, the length of the vertical component $PM$ is 1 and the length of the horizontal component $OP$ is $\frac{1}{cos\phi}$ on the projection, because the horizontal scale factor is $\frac{1}{\cos\phi}$ instead of 1. Next, consider $|PM'|=a$ where point $M'$ is the position of point $M$ on the projection and $0\leq{a}\leq1$. It follows that $|M'M|=1-a$.

Using figure 12, we aim to utilize all of its valuable information to calculate the maximum value of $\omega$. First, we notice that $
\tan\omega=\tan(\theta-\theta')$. Using the compound angle identity for tangent, we get 

$$\tan(\theta-\theta')=\frac{\tan\theta-\tan\theta'}{1+\tan\theta\tan\theta'}.$$

As can be seen from figure 12, $\tan\theta=\dfrac{1}{\frac{1}{\cos\phi}}=\cos\phi$ and $\tan\theta'=\dfrac{a}{\frac{1}{\cos\phi}}=a\cos\phi.$ Therefore, 

$$\frac{\tan\theta-\tan\theta'}{1+\tan\theta\tan\theta'}=\frac{\cos\phi-a\cos\phi}{1+\cos\phi\cdot{a}\cos\phi}=\frac{\cos\phi(1-a)}{1+a\cos{^2}\phi}.$$

Keep in mind that we are most interested in the maximum angle distortion. Thus, it would be plausible to differentiate the expression above with respect to $\phi$; the maximum may appear when the zeros of the derivative are retrieved.

$$\frac{d}{d\phi}\frac{\cos\phi(1-a)}{1+a\cos{^2}\phi}=(1-a)\frac{-\sin\phi\cdot(1+a\cos{^2}\phi)-(-\sin\phi\cos\phi-\sin\phi\cos\phi)\cdot\cos\phi\cdot{a}}{(1+a\cos{^2}\phi)^2}$$

$$=(1-a)\frac{-\sin\phi\cdot(1+a\cos{^2}\phi)+\sin{2\phi}\cdot\cos\phi\cdot{a}}{(1+a\cos{^2}\phi)^2}$$

$$=(1-a)\frac{-\sin\phi-a\sin\phi\cos{^2}\phi+a\sin{2\phi}\cos\phi}{(1+a\cos{^2}\phi)^2}$$

$$=(1-a)\frac{-\sin\phi-a\sin\phi\cos{^2}\phi+2a\sin\phi\cos\phi\cdot\cos\phi}{(1+a\cos{^2}\phi)^2}$$

$$=(1-a)\frac{2a\sin\phi\cos{^2}\phi-\sin\phi-a\cos{^2}\phi\sin\phi}{(1+a\cos{^2}\phi)^2}$$

$$=(1-a)\frac{(2a-a)\sin\phi\cos\phi-\sin\phi}{(1+a\cos{^2}\phi)^2}=(1-a)\frac{a\sin\phi\cos\phi-\sin\phi}{(1+a\cos{^2}\phi)^2}$$

The zeros of the derivative are the zeros of the numerator. Notice that $1-a$ is a constant since the expression $\dfrac{\cos\phi(1-a)}{1+a\cos{^2}\phi}$ is differentiated in respect to $\phi$.

$$a\sin\phi\cos\phi-\sin\phi=0\Leftrightarrow\sin\phi=a\sin\phi\cos\phi$$

Two solutions, $\sin\phi=0$ and $a\cos\phi=1$, are found. Solving $\sin\phi=0$ yields $\phi=0$ or $\phi=\pi$, but comparing the solutions to figure \ref{maximumomega}, it seems both are invalid: Since $|MP|\neq0$, $\phi\neq0$. Additionally, $\angle{MOP}\neq{\pi}$ since $\Delta{MOP}$ is a right traingle. Therefore, $\sin\phi\neq0$ and its solutions are invalid.

On the other hand, the expression $a\cos\phi=1$ illustrates a more interesting result. Solving the equation gives $\cos\phi=\frac{1}{a}$. Then, we can also substitute this result back to the original formula:

$$\frac{\cos\phi(1-a)}{1+a\cos{^2}\phi}=\frac{\frac{1}{a}(1-a)}{1+a(\frac{1}{a})^2}=\frac{\frac{1}{a}-1}{\frac{1}{a}+1}=\frac{1-a}{1+a}=\tan(\theta-\theta')=\tan\omega.$$

Solving the expression for $\omega$ results

$$\omega=\arctan\frac{1-a}{1+a}.$$

which shows the identical as substituting $\phi=0$ equation. Unlike the results found from substituting $\phi=0$, the result retrieved using $\cos\phi=\frac{1}{a}$ is valid, since $\cos\phi\neq0$ when $0\leq{a}\leq1$. 

\begin{figure}
\includegraphics[scale=0.5]{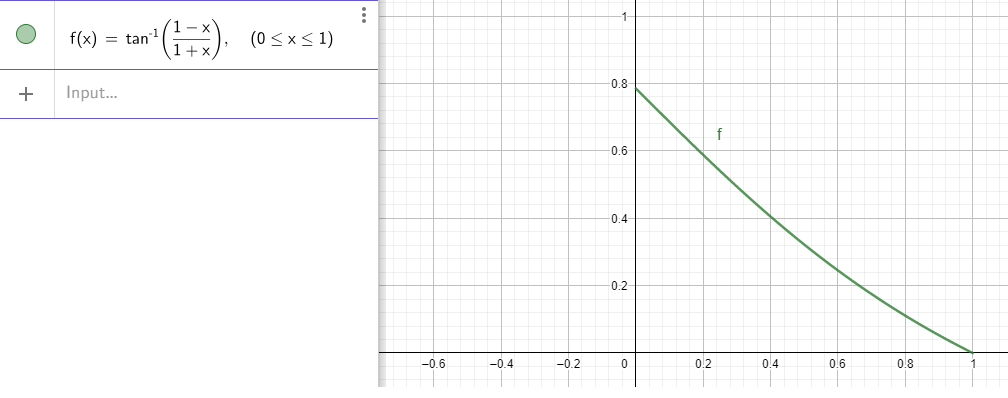}
\caption{The graph of $\arctan\frac{1-a}{1+a}$ in the range of $0\leq{a}\leq1$. It can be seen that $\arctan\frac{1-a}{1+a}$ is monotonously decreasing in the interval in question.}
\label{arctanx}
\end{figure}
By graphically investigating the expression $\arctan\dfrac{1-a}{1+a}$ in the range of $0\leq{a}\leq1$, we notice that 

$$a\in[0,1]\Rightarrow\arctan\frac{1-a}{1+a}\in[0,0.785...].$$

The range tells that the maximum angular distortion is $\omega\approx0.785$, which is obtained when $a=0$. This result implies that the distortion is at its minimum along the standard parallel of the projection, which is selected to be the equator.
 
Furthermore, we want to calculate $\dfrac{d}{da}\arctan\dfrac{1-a}{1+a}$, since the maximum of $\arctan\dfrac{1-a}{1+a}$ may occur at the zero of its derivative within $a\in[0,1]$. Notice that $\arctan\dfrac{1-a}{1+a}$ is in the form $f(g(a))$ where $f(a)=\arctan{a}$ and $g(a)=\dfrac{1-a}{1+a}$. Additionally, recall that $\dfrac{d}{da} f(g(a))=f'(g(a))+g'(a)$. Hence,

$$\frac{d}{da}\arctan\frac{1-a}{1+a}=\frac{\frac{-1-a-1+a}{(a+1)^2}}{\frac{(1-a)^2}{(1+a)^2}+1}=\frac{\frac{-2}{(a+1)^2}}{\frac{(1-a)^2+(1+a)^2}{(1+a)^2}}$$

$$=\frac{\frac{-2}{(a+1)^2}}{\frac{1-2a+a^2+1+2a+a^2}{(a+1)^2}}=\frac{-2}{2+2a^2}=\frac{-1}{1+a^2}.$$

It seems that $\dfrac{-1}{1+a^2}$ has no real zeroes in its defined domain ($a\in R$). Hence, the expression has no zeroes within the range $a\in [0,1]$. Moreover, $\dfrac{-1}{1+a^2}<0$ everywhere when $a\in[0,1]$; Thus, $\dfrac{-2}{(1+a)^2}$ is decreasing when the value of $a$ increases, therefore the value $a=0$ will give the maximum angular distortion. In conclusion, $\omega=0.785$ is indeed, the maximum angular distortion possible. 

From page 12, we recall that the rate of horizontal distortion will increase without bound at polar areas; this phenomenon can also be shown using the concept of Tissot's indicatrix. By generalizing the lengths of $OP$ and $PM$ in figure 12 for all possible longitudinal and latitudinal coordinates, we can easily calculate the length of $OM$, which is the hypotenuse of the triangle. Let $OP=\lambda\dfrac{1}{\cos\phi}$ and $MP=\phi$, where $\lambda$ and $\phi$ are the differences in longitudinal/latitudinal coordinates and $-\pi\leq\lambda\leq\pi$ and $-\frac{\pi}{2}\leq\phi\leq\frac{\pi}{2}$. According to Pythagoras' rule, $|OM|=\sqrt{\phi^2+\lambda{^2}\sec{^2}\phi}$.

However, the distance $OM$ is distorted on the map projection: The scale factors are not taken into account in our model (page \pageref{Front view of 3D figure}, figure \ref{Front view of 3D figure}). Therefore, the projection of $OP$- which is $OP'$- is equal to $\lambda$. Subsequently, the projected line $OM'$ will be equal to $\sqrt{\lambda^2+\phi^2}$ in the 2D plane.

As mentioned in page 18, the maximum distortion in distance is equal to $||OM|-|OM'||=\sqrt{\phi^2+\lambda{^2}\sec{^2}\phi}-\sqrt{\lambda^2+\phi^2}.$ Examining the expression within the range $-\frac{\pi}{2}\leq\phi\leq\frac{\pi}{2}$, we notice that when $\phi=0$, $\sqrt{\phi^2+\lambda{^2}\sec{^2}\phi}-\sqrt{\lambda^2+\phi^2}=0$; the distance on the equator - our standard parallel - is perfectly preserved. Additionally, when $\phi\rightarrow
\pm\frac{\pi}{2}$, $\sqrt{\phi^2+\sec{^2}\phi}-\sqrt{1+\phi^2}\rightarrow\infty$, implying that the distortion in distance is infinitely large in polar areas. 

In conclusion, this section proves that flexion arises due to a doubtful assumption: The projection assumes all parallels to be homogeneous straight lines equal in length. Unlike meridians, not all parallels are equal in length, as figure \ref{Standard Parallels = circles} demonstrates. To make all parallels equal in length, points on the Earth are projected onto the surface of the cylinder as figure \ref{Front view of 3D figure} illustrates. Thus, as proven on page 10, horizontal distances are distorted by a factor of $\dfrac{1}{\cos\phi}$. This operation causes artificial bending, due to which the shape of the indicatrices change. As shown in pages \pageref{indicatrix}-\pageref{Extramashup}, the bending of Tissot's indicatrices will determine the extent of distortions in distance (values of $\sqrt{\phi^2+\sec{\phi^2}}-\sqrt{1+\phi^2}$ when $-\frac{\pi}{2}\leq\phi\leq\frac{\pi}{2}$). Moreover, the errors in distance are dependent only on one variable - $\phi$, and cannot be alleviated by changing other variables. Therefore, as the distance distortion in the horizontal axis exists and cannot be canceled out by other bending, distance cannot be completely preserved.

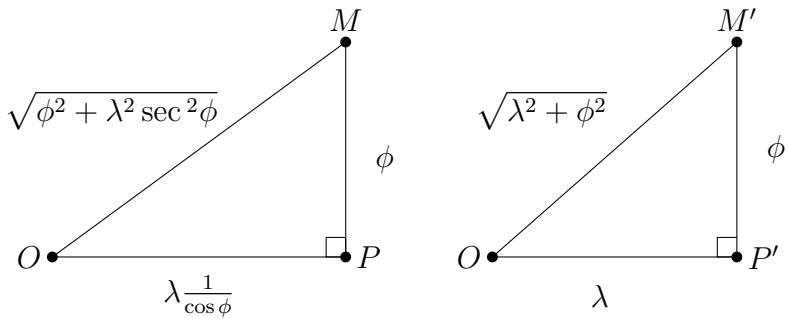
\begin{figure}
\begin{tikzpicture}[x=0.13cm,y=0.13cm]

\draw (-10,0)--(20,0);
\draw (20,0)--(20,22);
\draw (20,22)--(-10,0);
\draw (20,0)--(18,0)--(18,2)--(20,2)--(20,0);
\draw [fill] (-10,0) circle [radius=0.5] node[left] {$O$};
\draw [fill] (20,22) circle [radius=0.5] node[above] {$M$};
\draw [fill] (20,0) circle [radius=0.5] node[right] {$P$};
\node at (5,-4) {$\lambda\frac{1}{\cos\phi}$};
\node at (24,10) {$\phi$};
\node at (-4,15) {$\sqrt{\phi^2+\lambda{^2}\sec{^2}\phi}$};

\draw (35,0)--(60,0);
\draw (60,0)--(60,22);
\draw (60,22)--(35,0);
\draw (60,0)--(58,0)--(58,2)--(60,2);
\draw [fill] (35,0) circle [radius=0.5] node [left] {$O$};
\draw [fill] (60,0) circle [radius=0.5] node [right] {$P'$}; 
\draw [fill] (60,22) circle [radius=0.5] node [above] {$M'$};

\node at (46,-4) {$\lambda$};
\node at (64,11) {$\phi$};
\node at (40,15) {$\sqrt{\lambda^2+\phi^2}$};
\end{tikzpicture}
\caption{A generalized illustration of the lengths of the sides $OP$, $MP$ and $OM$ in figure 12. On the projection, the projected side $OM'$ forms another relationship with sides $OP'$ and $M'P'$, from which the distortion in distance can be demonstrated.} 
\end{figure}

\newpage
\subsection*{Trigonometrical analysis}

In this section, we will further examine the relationship between the "straight line" distance between two points on an equidistant cylindrical projection and the real distance on the Earth's surface predominately using trigonometry and sinusoidal identities. 

\begin{figure}

\tdplotsetmaincoords{70}{110}
\begin{tikzpicture}[tdplot_main_coords]
\draw [->] (0,0,0)--(5,0,0) node[anchor=north west]{$x$} ;
\draw [->] (0,0,0)--(0,5,0) node[anchor=north east]{$y$} ;
\draw [->] (0,0,0)--(0,0,5) node[anchor=south] {$z$};
\draw [dashed] (0,0,0) circle [radius=1];
\draw [dashed] (0,0,0.85) to [out=90,in=90] (0,-1,0);
\draw [dashed] (0,0,1) to [out=0,in=90] (0,1,0);
\draw [dashed] (0,1,0) to [out=-90,in=0] (0,0,-1);
\draw [dashed] (0,0,-1) to [out=-180,in=-90] (0,-1.05,0);
\draw [green](0,0,0) -- (0,0.707,0.707);
\draw [fill] (0,-0.3,0.5) circle [radius=0.1];
\draw [fill] (0.6,0.2,-0.1) circle [radius=0.1];
\node at (2,-0.5,0) {$1$};
\node at (0,1.5,0.2) {$1$};
\node at (0,1.414,1.414) {$r=1$};
\node at (-2,-1.5,0.1) {P$_1$};
\node at (2,0.5,0) {P$_2$};

\end{tikzpicture}
\caption{a unit sphere in xyz-coordinate system with arbitrary points P$_1$ and P$_2$ on its surface. It can be understood that the shortest distancefrom P$_1$ to P$_2$ is a curvature along the surface of the sphere.}
\label{unit sphere}
\end{figure}
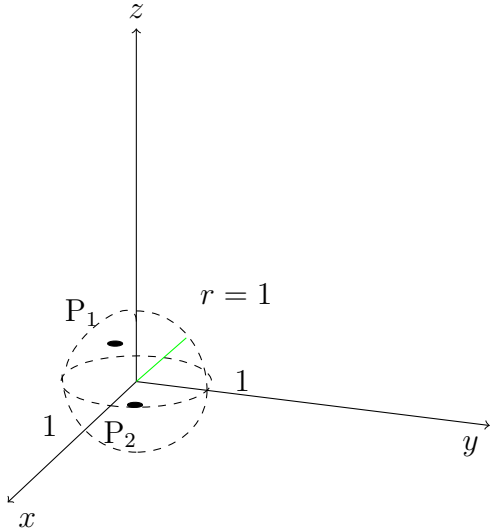

First, let sphere with a radius of 1 be in a xyz-coordinate system, as shown in figure \ref{unit sphere}. Arbitrarily selected points P$_1$ and P$_2$ lie on the surface of the sphere. If an equidistant cylindrical projection is generated from the sphere, the distance between P$_1$ and P$_2$, would be $\sqrt{(\lambda_1-\lambda_2)^2+(\phi_1-\phi_2)^2}$, according to the Pythagoras' formula. (See figure \ref{P1 and P2})

However, we can observe that the \textit{de facto} shortest path connecting P$_1$ and P$_2$ is not a line but a curve along the surface of the sphere, since it is impossible to travel through the Earth. Therefore, we cannot purely examine the straight-line Euclidean distance which we calculate using the formulas of equidistant cylindrical projection. We can convert the Cartesian coordinates into \textbf{spherical coordinates}\footnote{Paul Dawkins, \textsl{Calculus III: Spherical coordinates}, retrieved June 2019. \\
\url{http://tutorial.math.lamar.edu/Classes/CalcIII/SphericalCoords.aspx}}, which can be used to treat the route taken from P$_1$ to P$_2$ as a curve. Using the spherical coordinate system, we want to  compare the calculated "straight-line" distance between two points on equidistant cylindrical projection to the real distance between those two points on the surface of the Earth. 

\subsubsection*{Spherical coordinates}

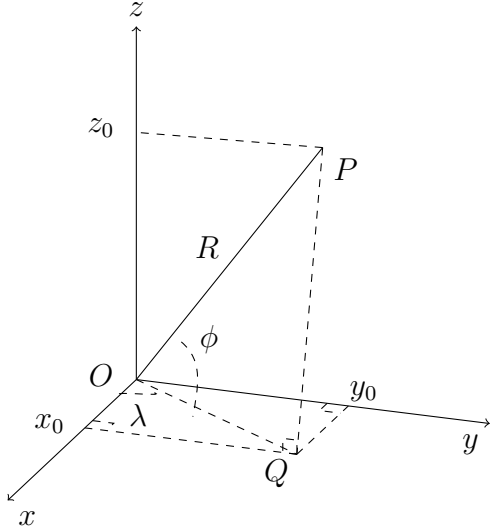
\begin{figure}

\tdplotsetmaincoords{70}{110}
\begin{tikzpicture}[tdplot_main_coords]
\draw [->] (0,0,0)--(5,0,0) node[anchor=north west]{$x$} ;
\draw [->] (0,0,0)--(0,5,0) node[anchor=north east]{$y$} ;
\draw [->] (0,0,0)--(0,0,5) node[anchor=south] {$z$};
\draw (0,0,0)--(1,3,4) node[anchor=north west] {$P$};
\draw [dashed] (1,3,4)--(2,3,0);
\draw [dashed] (0,0,0)--(2,3,0);
\draw [dashed] (2,3,0)--(0,3,0);
\draw [dashed] (2,3,0)--(2,0,0);
\draw [dashed] (1,3,4)--(0,0,3.5);
\draw [dashed] (2,3,0)--(2,2.8,0)--(2,2.8,0.2)--(2,3,0.2)--(2,3,0);
\draw [dashed] (0.4,0.4,0) to [out=-20,in=150] (0.7,0,0);
\draw [dashed] (1,1,1) to [out=-40,in=80] (1.1,1.2,0);
\draw [dashed] (1.7,0,0)--(1.7,0.3,0)--(2,0.3,0);
\draw [dashed] (0,2.7,0)--(0.3,2.7,0)--(0.3,3,0);
\node at (1.3,0.5,0) {$\lambda$};
\node at (1.3,1.5,1.2) {$\phi$};
\node at (0,1,2) {$R$};
\node at (0,-0.5,0) {$O$};
\node at (2,-0.5,0) {$x_0$};
\node at (-0.6,3,0) {$y_0$};
\node at (0,-0.5,3.5) {$z_0$};
\node at (2.8,3,0) {$Q$};

\end{tikzpicture}
\caption{An illustration of spherical coordinates:}
\label{sphericalcoordinates}
\end{figure}

We can convert \footnote{University of Utah, \textsl{Cylindrical and spherical coordinates}, retrieved June 2019.

\url{https://www.math.utah.edu/lectures/math2210/9PostNotes.pdf}} the Cartesian coordinates $x$, $y$ and $z$ into spherical coordinate notation. 

In figure \ref{sphericalcoordinates}, let $\vec{OP}$ be a vector in 3D plane and the Euclidean distance $|\vec{OP}|=R.$ It can be interpreted that the $x$, $y$ and $z$ components of $\vec{OP}$ are $\vec{Ox_0}$, $\vec{Oy_0}$ and $\vec{QP}$ respectively. Let the $x$, $y$ and $z$ coordinates of point $P$ and $Q$ be ($x_0$, $y_0$, $z_0$) and ($x_0,y_0,0$) respectively. The spherical coordinate system is able to represent those Cartesian coordinates in terms of angles $\lambda$ and $\phi$. As figure 15 states, $\Delta{OPQ}$, $\Delta{Ox{_0}Q}$ and $\Delta{Oy{_0}Q}$ are all right triangles. Hence, trigonometrical identities can be used to solve for the magnitudes of the components of vector $\vec{OP}$, which are $\vec{Ox_0}$, $\vec{Oy_0}$ and $\vec{PQ}$.

To solve for $|\vec{Ox_0}|$, we notice that $|\vec{Ox_0}|=\cos\lambda|\vec{OQ}|$, where $|\vec{OQ}|=\cos\phi|\vec{OP}|$. By substituting $|\vec{OP}|=R$ and simplifying the expression, 

$$|\vec{Ox_0}|=\cos\lambda|\vec{OQ}|=\cos\lambda\cos\phi|\vec{OP}|=R\cos\lambda\cos\phi.$$

Since both $\Delta{Ox{_0}Q}$ and $\Delta{Oy{_0}Q}$ are right triangles, share a common side $OQ$, and sides $Oy_0\|Qx_0$ and $Ox_0\|y{_0}Q$, it seems that $\Delta{Ox{_0}Q}$ and ${Oy{_0}Q}$ are similar triangles. Hence, because of alternating angles, $\angle{OQy_0}=\lambda$. Thus, the other components, $|\vec{Oy_0}|$ and $|\vec{PQ}|$, can be calculated:

$$|\vec{Oy_0}|=\sin\lambda|\vec{OQ}|=R\sin\lambda\cos\phi$$

$$|\vec{PQ}|=R\sin\phi$$

In conclusion, it seems that the spherical coordinates are for all points with Cartesian coordinates $(x,y,z)$ can be written as

$$x=R\cos\lambda\cos\phi$$

$$y=R\sin\lambda\cos\phi$$

$$z=R\sin\phi$$

Using the concept of spherical coordinates, the distance between points P$_1$ and P$_2$ can be calculated using the longitude $\lambda$ and latitude $\phi$. First, we want to calculate the straight-line distance between P$_1$ and P$_2$ (Cf. figure \ref{P1 and P2}) by treating the two points as projections onto a unit sphere. Using the Pythagorean theorem, the straight-line distance between P$_1$ and P$_2$ in the unit sphere can be deduced: \footnote{Kansas State University, \textsl{Distance between points on the Earth's surface}, retrieved June 2019.

\url{https://www.math.ksu.edu/~dbski/writings/haversine.pdf}} $$d^2=(x_1-x_2)^2+(y_1-y_2)^2+(z_1-z_2)^2.$$ By converging the variables $x$, $y$ and $z$ into spherical coordinates we get $$d^2=(x_1-x_2)^2+(y_1-y_2)^2+(z_1-z_2)^2$$
$$=R^2[(\cos\lambda_1\cos\phi-\cos\lambda_2\cos\phi_2)^2+(\sin\lambda_1\cos\phi_1-\sin\lambda_2\cos\phi_2)^2+(\sin\phi_1-\text{sin}\phi_2)^2]$$
$$=R^2[\cos^2\lambda_1\cos^2\phi_1-2\cos\lambda_1\cos\lambda_2\cos\phi_1\cos\phi_2$$ $$+\cos^2\lambda_2\cos^2\phi_2+\sin^2\lambda_1\cos^2\phi_1-2\sin\lambda_1\sin\lambda_2\cos\phi_1\cos\phi_2$$ $$+\sin^2\lambda_2\cos^2\phi_2+\sin^2\phi_1-2\sin\phi_1\sin\phi_2+\sin^2\phi_2.]$$

Recall that $\sin^2{a}+\cos^2{a}=1$ and $\cos(a-b)=\cos{a}\cos{b}+\sin{a}\sin{b}$. Thus, rearranging and simplifying the terms results us 

$$d^2=R^2[\cos^2\lambda_1\cos^2\phi_1+\sin^2\lambda_1\cos^2\phi_1+\cos^2\lambda_2\cos^2\phi_2+$$ $$\sin^2\lambda_2\cos^2\phi_2+\sin^2\phi_1+\sin^2\phi_2-2\cos\lambda_1\cos\lambda_2\cos\phi_1\cos\phi_2-$$ $$2\sin\lambda_1\sin\lambda_2\cos\phi_1\cos\phi_2-2\sin\phi_1\sin\phi_2]$$
$$=R^2[\cos^2\phi_1(\cos^2\lambda_1+\sin^2\lambda_2)+\cos^2\phi_2(\cos^2\lambda_2+\sin^2\lambda_2)+   \sin^2\phi_1+\sin^2\phi_2$$
$$-2\cos\phi_1\phi_2cos(\lambda_1-\lambda_2)-2\sin\phi_1\sin\phi_2]$$
$$=R^2[\cos^2\phi_1+\sin^2\phi_1+\cos^2\phi_2+\sin^2\phi_2-2\cos\phi_1\phi_2\cos(\lambda_1-\lambda_2)-2\sin\phi_1\sin\phi_2]$$
$$=R^2(2-2\cos\phi_1\cos\phi_2\cos(\lambda_1-\lambda_2)-2\sin\phi_1\sin\phi_2).$$

Therefore, by examining the sphere in spherical coordinates, it is determined that the straight-line distance between two points on the surface of a sphere is

$$d=\sqrt{R^2(2-2\cos\phi_1\cos\phi_2\cos(\lambda_1-\lambda_2)-2\sin\phi_1\sin\phi_2)}$$ $$=R\sqrt{2-2\cos\phi_1\cos\phi_2\cos(\lambda_1-\lambda_2)-2\sin\phi_1\sin\phi_2},$$ where $R$ is the radius of the Earth.

\subsubsection*{Relationship between straight-line distance and real distance}

\begin{figure}
\begin{tikzpicture}
\draw [-] (0,0) -- (3,1) node[above] {P$_1$};
\draw [-] (0,0) -- (3,-1) node[below] {P$_2$};
\draw [dashed] (3,1) -- (3,-1);
\draw [dashed] (0,0) -- (3,0) node[right] {$d$};
\draw [thick] (3,1) to [out=0,in=0] (3,-1) node at (4,0) {$D$};
\draw [thick] (1,0.33) to [out=0,in=0] (1,-0.33) node at (-0.2,0) {$\alpha$};
\draw [thick] (0.9,0.28) to [out=0,in=0] (1,0) node at (0.6,0.7) {$\frac{\alpha}{2}$};
\draw [dashed] (3,0) -- (2.8,0) -- (2.8,0.2) -- (3,0.2);
\node at (1.5,0.85) {$R$};
\node at (1.5,-0.85) {$R$};
\node at (2.7,0.5) {$\frac{d}{2}$};

\end{tikzpicture}
\caption{A sketch of points P$_1$ and P$_2$, their "straight-line" distance $d$ and their real distance $D$.}
\label{fig:16}
\end{figure}
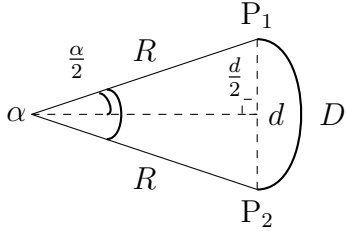

Next, we want to examine the relationship between the real surface distance between two points on a sphere and the straight-line distance that we obtained (see \texttt{spherical coordinates}). In figure~\ref{fig:16}, we form a model containing both the straight-line distance denoted as $d$, and the real distance, denoted as $D$, on the surface of the sphere between points P$_1$ and P$_2$. The radius of the sphere is assumed to be $R$. The angle bisector of $\alpha$ is perpendicular to line $d$ and bisects $d$, because in the triangle  $\Delta$P$_1$P$_2$ the angle $\angle$ $\alpha$ is isosceles.

From the figure, we observe that the arc length - which is the real distance between P$_1$ and P$_2$ is $D$, is equal to $R\alpha$. We may link $d$ and $D$ together if we are able to calculate $\alpha$. Using trigonometric identities, we are able to write $\alpha$ with respect to $d$, the "straight-line" distance between P$_1$ and P$_2$. First, we see that $\sin\frac{\alpha}{2} = \frac{d}{2R}$. By using double-angle identities, $$\sin\alpha = \sin (2\cdot\frac{\alpha}{2}) = 2\sin\frac{\alpha}{2}\cos\frac{\alpha}{2}$$

Because $2\sin\frac{\alpha}{2} = 2\cdot\frac{d}{2R} = \frac{d}{R}$ and $\cos\frac{\alpha}{2}=\sqrt{1-\sin^{2}\frac{\alpha}{2}}=\sqrt{1-(\frac{d}{2R}})^2$, we can simplify $\sin\alpha$ as $$\sin\alpha=\frac{d}{R}\cdot\sqrt{1-(\frac{d}{2R})^2}=\frac{d}{R}\sqrt{1-\frac{d^2}{4R^2}}=\frac{d}{R}\sqrt{\frac{4R^2-d^2}{4R^2}}=\frac{d}{2R^2}\sqrt{4R^2-d^2}.$$

Thus, as $D=R\alpha,$ we conclude that $$D=R\alpha=R\arcsin(\frac{d}{2R^2}\sqrt{4R^2-d^2}).$$

Hence, there exists a relationship between the "straight-line" Euclidean distance between two arbitrary points and the real distance between those two points along the surface of the Earth. 
We can now use this relationship to calculate the real distance between the points P$_1$ = (24.3$^{\circ}$E, 23.4$^{\circ}$N) $\approx$ (0.424,0.408) and P$_2$ = (39.2$^{\circ}$W, 3.67$^{\circ}$S) $\approx(-0.492,0.674)$ on page 11. Recall that the straight-line distance between P$_1$ and P$_2$ is calculated using the expression $d=R(2-2\cos\phi_1\cos\phi_2\cos(\lambda_1-\lambda_2)-2\sin\phi_1\sin\phi_2)$. Thus, 
\newpage
$$d_{P_1,P_2}=R\sqrt{2-2\cos\phi_1\cos\phi_2\cos(\lambda_1-\lambda_2)-2\sin\phi_1\sin\phi_2}$$ $$=6371\cdot\sqrt{(2-2\cos0.408\cos(-0.064)\cos(0.424-(-0.687))-2\sin0.408\sin(-0.064)}$$ $$\approx 6371\cdot\sqrt{2-0.811-0.050}\approx 7091.584111$$

And thus, $$D=R\arcsin(\frac{d}{2R^2}\sqrt{4R^2-d^2})$$ $$=6371\cdot\arcsin(\frac{7091.584111}{2\cdot6371^2}\sqrt{4\cdot6371^2-7091.584111^2})
\approx 6371\cdot1.180\approx7520$$

Compare the result above with the actual distance\footnote{Movable type scripts, \textsl{Calculate distance, bearing and more between latitude/longitude points}

https://www.movable-type.co.uk/scripts/latlong.html} in page 15 ($d=7502$) gives us a percentage error of $\frac{7520-7502}{7502}\cdot100\%=0.240\%.$ 

The percentage error in distance calculated using equidistant cylindrical projection formulas is $2.28\%$ (see page 10), which is significantly greater than using the formula $D=R\arcsin(\frac{d}{2R^2}\sqrt{4R^2-d^2})$ to calculate the same distance. Looking back to the assumptions made when creating the model, it seems that a key feature of the projection is based on a fallacy: the straight-line distance between two arbitrary points on the projection (figure \ref{P1 and P2}) is assumed to always correspond to the real distance between those two points on Earth. There are many examples that disprove the assumption, which as shown below.

As demonstrated by the calculations on page 13, the real distance between P$_1$ and P$_2$ is shorter than what the projection claims. Now, we are going to examine two more pairs of points that further demonstrates that the projection is not be able to present the shortest possible distance between two points.

\textbf{Example 1.} Consider arbitrarily selected points A $=(41.0^{\circ}$W, $9.2^{\circ}$S) $\approx (-0.716,-0.161)$ and B $=(48.1^{\circ}$E, $21.8^{\circ}$N) $\approx(0.840,0.380)$. According to equidistant cylindrical projection, $$d_{A,B}=\sqrt{(R\lambda_1-R\lambda_2)^2+(R\phi_1-R\phi_2)^2}=10495.$$

However, by using the formula calculated using spherical coordinates, the "straight-line distance" between points A and B are found to be $$d=R\sqrt{2-2\cos\phi_1\cos\phi_2\cos(\lambda_1-\lambda_2)-2\sin\phi_1\sin\phi_2}=9214$$

and the real distance separating A and B along the Earth is $$D=R\arcsin(\frac{d}{2R^2}\sqrt{4R^2-d^2})=9715<10495.$$

\textbf{Example 2.} Consider another set of randomly selected points C = $(170.8^{\circ}$E, $65.1^{\circ}$N) = (2.98,1.13) and D = $(152.7^{\circ}$W, $64.7^{\circ}$N) = (-2.67,1.13). According to equidistant cylindrical projection, $$d_{C,D}=\sqrt{(R\lambda_1-R\lambda_2)^2+(R\phi_1-R\phi_2)^2}=35760$$

Once again, by examining the distance using the formula deduced using spherical coordinates, we get $$d=R\sqrt{2-2\cos\phi_1\cos\phi_2\cos(\lambda_1-\lambda_2)-2\sin\phi_1\sin\phi_2}=1689.$$

and the real distance between points C and D is $$D=R\arcsin(\frac{d}{2R^2}\sqrt{4R^2-d^2})=1694<<35760.$$

From the two examples above, we can deduce that equidistant cylindrical projection does not have the ability to perfectly represent the actual distance between any two points on the surface of the Earth. The projection is able to represent a straight-line route connecting two points as seen in figures 5 and 6, but the length of the route may deviate from what the actual shortest route is. Therefore, it can be said that the formula $D=R\arcsin(\frac{d}{2R^2}\sqrt{4R^2-d^2})$ is able to generate a more accurate result for the real distance between P$_1$ and P$_2$ than the formulas of the equidistant cylindrical projections. Hence, it can be concluded that equidistant cylindrical projections fails to demonstrate the actual route between two arbitrary points on the surface of the Earth as curves, thus distances become less accurately preserved in equidistant cylindrical projections.

In conclusion, the straight-line distance which equidistant cylindrical projection illustrates differs from the real distance on the surface of the sphere; as a result, we can deduce that equidistant cylindrical projection is unable to perfectly represent the geodesic distance between two points: the projection must contain distortion in distance. 

\newpage
\section*{Conclusion}
From this research, we can conclude that the equidistant cylindrical projection cannot fully preserve distance between two arbitrary points on the surface of the Earth. We have successfully found out how can distortions in distance be detected, and to what extent does the projected distance vary from the real distance. Finally, this research demonstrates that the distortion in distance is minimum at the equator and the maximum at polar areas in equidistant cylindrical projection. Thus, the aim of this work is achieved.

It is undeniable that the projection has some distance-preserving qualities. However, the design of the projection has major flows that cause distortion in distance. By examining the distance between two arbitrary points on the modeled projection using different methods, such as the Tissot's indicatrix and the spherical coordinates, it is deduced that the projection contains a relevant amount of distortion in distance. 

In sections \textsl{Geometrical analysis} and \textsl{Trigonometrical analysis}, the degree of distance distortion in the equidistant cylindrical projection is shown from different aspects. Based on the result of this research, it can be said that there are two main questionable assumptions made during the map-making phase:

\begin{enumerate}

\item{All parallels are assumed to have the same length.}

\item{The "straight-line" route between two points is assumed to be the shortest possible route between those two points on the surface of the Earth.}
\end{enumerate}

As explored in this research, those two assumptions are directly responsible for the inaccurate illustration of distance between two points on the projection. Additionally, there are other assumptions made in the making of equidistant cylindrical projections, such that the Earth is assumed to be perfectly spherical and the radius of the Earth is exactly 6371 kilometers. Further investigations can be conducted to determine the effect of different assumptions on the distance-preserving ability of equidistant cylindrical projections.
\newpage
\begin{center}
\section*{Appendix}

\end{center}

\subsection*{Glossary}
\begin{itemize}

\item{Map projection: A representation of the properties of the spherical Earth on a flat surface.}

\item{Cylindrical projection: A projection that is formed by centering a sphere into a right cylinder and individual points on the sphere are transferred onto the side surface of the cylinder. The end product is a rectangular sheet, parallels and meridians intersect each other in right angles.}

\item{Equidistant cylindrical projection: A type of cylindrical projection that is meant to describe distances between any two arbitrary points, which should correspond to their actual distance on the Earth.}

\item{Meridian: Any great circle of the earth that vertically passes through the poles and any given point on the earth's surface.}

\item{Standard parallel: Any great circle of the earth that is horizontally parallel to the equator.}

\item{Scale factor: A value that corresponds to the distance of a specific route taken to connect two points on a map projection in respect of the real distance of those two points on the Earth.}

\item{Distortion: Any angular, areal or distant error that is present in map projections.}

\item{Geodesic: Denotation of the shortest possible distance between two points on a sphere or other curved surface.}

\item{Tissot's indicatrix: An infinitely small circle drawn on the surface of Earth. The circle's shape is distorted once the Earth is projected onto a flat sheet.}

\item{Flexion: Artificial bending of structures.}

\end{itemize}

\newpage

\subsection*{General model for equidistant cylindrical projection}
\textsl{Note:} the formulas listed in page 10 are those that are used in this research, but they are not the general formulas for equidistant cylindrical projections; the central meridian and standard parallel are selected to be specific values. According to \footnote{John P. Snyder (1987), \textsl{Map projection - A working manual}}Snyder (1987), the generalized formulas for any equidistant cylindrical projections are

$$y=R\phi$$ 

$$x=R(\lambda-\lambda_{0}){\text{cos}\phi_{1}}$$

$$h=1$$

$$k=\dfrac{\text{cos}\phi_{1}}{\text{cos}\phi}$$

$$-\frac{\pi}{2}\leq\phi\leq\frac{\pi}{2}, -\pi\leq\lambda\leq\pi$$

Where 

\begin{itemize}

\item{$\lambda_0$ is the selected central meridian}

\item{$\phi_1$ is the selected standard parallel.}
 
\end{itemize} 
 
The other variables are explained in page 11. 
 
In Snyder's model, the values of $k$ - which causes distortion in distance in the first place - and $x$ - which, due to $k$, suffers from artificial bending - are written to be dependent on the central meridian $\lambda_0$ and standard parallel $\phi_1$, which are arbitrarily selected.

\newpage
\subsection*{The Haversine formula}

In an earlier subsection, the shortest route between two points on the surface of the sphere is examined as a curve along the surface of the Earth. This distance can also be denoted\footnote{Kansas State University, \textsl{Distance between points on the Earth's surface}, retrieved June 2019.

\url{https://www.math.ksu.edu/~dbski/writings/haversine.pdf}} using the \texttt{Haversine formula}, which is defined as $$\text{haver}\sin{\alpha}=\sin^2\frac{\alpha}{2}=\frac{1-\cos\alpha}{2},$$ from which we determine $$\cos\alpha=1- 2 \text{haver}\sin\alpha,$$ where $\alpha$ is the angle in figure \ref{fig:16}.

In this section, the real distance between two arbitrary points on the surface of the Earth will be further examined using the Haversine formula, which is already briefly discussed in an earlier section. Using the expression $\cos{\alpha}=1-2\text{haver}\sin{\alpha}$, we can convert the relationship of the "straight-line" distance between two points and the actual distance between them on the Earth.

From page 24, we recall that $d^2=R^2(2-2\cos\phi_1\cos\phi_2\cos(\lambda_1-\lambda_2)-2\sin\phi_1\sin\phi_2).$ Substituting $\cos{\lambda_{1}-\lambda_{2}}=1-2\text{haver}\sin{\lambda_{1}-\lambda_{2}}$, the expression can be rewritten as

$$d^2=R^2(2-2\cos\phi_1\cos\phi_2\cos(\lambda_1-\lambda_2)-2\sin\phi_1\sin\phi_2)$$ $$=R^2(2-2\cos\phi_1\cos\phi_2[1-2\text{haver}\sin(\lambda_1-\lambda_2)]-2\sin\phi_1\sin\phi_2$$ $$=R^2(2-2\cos\phi_1\cos\phi_2-4\cos\phi_1\cos\phi_2\text{haver}\sin(\lambda_1-\lambda_2)-2\sin\phi_1\sin\phi_2)$$ $$=R^2(2-2\cos(\phi_1-\phi_2)+4\cos\phi_1\cos\phi_2\text{haver}\sin(\lambda_1-\lambda_2)$$ $$=R^2(4\cdot\frac{1-\cos(\phi_1-\phi_2)}{2}+4\cos\phi_1\cos\phi_2\text{haver}\sin(\lambda_1-\lambda_2)$$ $$=R^2(4\text{haver}\sin(\lambda_1-\lambda_2)+4\cos\phi
_1\cos\phi_2\text{haver}\sin(\lambda_1-\lambda_2).$$ \\

By dividing both sides of the equation by $4R^2$, it seems that  $$(\frac{d}{2R})^2=\text{haver}\sin\alpha=\text{haver}\sin(\phi_1-\phi_2)+\cos\phi
_1\cos\phi_2\text{haver}\sin(\lambda_1-\lambda_2).$$

Additionally, since $$d^2=4R^2\text{haver}\sin\alpha$$ and $$\frac{d}{2R^2}\sqrt{4R^2-d^2}=\sqrt{\frac{d^2}{4R^4}\cdot(4R^2-d^2)}=\sqrt{\frac{d^2}{R^2}-4\frac{d^4}{4R^4}},$$ 
it can be deduced that 
$$\frac{d}{2R^2}\sqrt{4R^2-d^2}=\sqrt{\frac{4R^2\text{haver}\sin\alpha}{R^2}-\frac{4\cdot{16}R^4\text{haver}\sin\alpha}{4R^4}}$$ $$=\sqrt{4\text{haver}\sin\alpha-16\text{haver}\sin^{2}\alpha}=2\sqrt{\text{haver}\sin\alpha-4\text{haver}\sin^{2}\alpha}.$$

The results retrieved above can be effectively utilized to digitally generate routes between two points on the surface of an object with eccentricity (sphere, ellipsoid etc.).   For example, some mapping algorithms, such as \textsl{Movable Type Script}\footnote{Movable type scripts, \textsl{Calculate distance, bearing and more between latitude/longitude points} 

https://www.movable-type.co.uk/scripts/latlong.html}, use the Haversine formula to accurately calculate the distance between any two points on the surface of the Earth using the shortest possible route between them. Therefore, the Haversine formula has potential of contributing tremendous real-life application value, hence revamping the techniques used in the map-making industry.
\newpage

\subsection*{Bibliography}

\begin{enumerate}
\subsubsection*{Books}

\item{Ivars Peterson, \textsl{The mathematical tourist: snapshots of modern mathematics}}

\subsubsection*{Theses}

\item{John P. Snyder (1987), \textsl{Map projections - A working manual}

\url{https://pubs.usgs.gov/pp/1395/report.pdf}} 

\item{David M. Goldberg, J. Richard Gott III (2007), \textsl{Flexion and Skewness in Map Projections of the Earth}

\url{https://www.physics.drexel.edu/~goldberg/projections/goldberg_gott.pdf}}

\item{Wenping Jiang, Jin Li (2014), \textsl{The Effects of Spatial Reference Systems on the Predictive Accuracy of Spatial Interpolation Methods}

\url{https://www.researchgate.net/publication/259848543_The_Effects_of_Spatial_Reference_Systems_on_the_Predictive_Accuracy_of_Spatial_Interpolation_Methods}}

\item{Bernhard Jenny, Bojan Savric (2017), \textsl{Combining world map projections}

\url{https://www.researchgate.net/publication/315860317_Combining_World_Map_Projections}}

\subsubsection*{Teaching and learning materials}

\item{Wikipedia, \textsl{Equidistant cylindrical projection,} retrieved December 2018.}

\item{Ball State University, \textsl{Maps and cartography; map projections}, retrieved February 2019.

\url{http://cardinalscholar.bsu.edu}}

\item{Benjamin H{\o}yer, \textsl{Cartography}, retrieved March 2019.

\url{http://www-history.mcs.st-and.ac.uk/Projects/Hoyer/S4.html}}

\item{James S. Aber, \textsl{Brief history of maps and cartography}, retrieved March 2019.

\url{http://academic.emporia.edu/aberjame/map/h_map/h_map.htm}}

\item{City University of New York, \textsl{Projection parameters}, retrieved March 2019.

\url{http://www.geography.hunter.cuny.edu/~jochen

/gtech361/lectures/lecture04/concepts/Map}}

\item{California State University Long Beach,\textsl{The geographic grid}, retrieved March 2019.

\url{https://web.csulb.edu/~rodrigue/geog140/lectures/geographicgrid.html}}

\item{GIS Geography, \textsl{Map distortions with Tissot's indicatrix}, retrieved April 2019.

\url{https://gisgeography.com/map-distortion-tissots-indicatrix/}}

\item{Rice University, \textsl{Mapping the sphere}, retrieved April 2019.

\url{https://math.rice.edu/~polking/cartography/cart.pdf}}

\item{Kansas State University, \textsl{Distance between points on the Earth's surface}, retrieved June 2019.

\url{https://www.math.ksu.edu/~dbski/writings/haversine.pdf}}

\item{Paul Dawkins, \textsl{Calculus III: Spherical coordinates}, retrieved June 2019.

\url{http://tutorial.math.lamar.edu/Classes/CalcIII/SphericalCoords.aspx}}

\item{University of Utah, \textsl{Cylindrical and spherical coordinates}, retrieved June 2019.

\url{https://www.math.utah.edu/lectures/math2210/9PostNotes.pdf}}

\item{Eotvos Lorand University, \textsl{Cylindrical projections}, retrieved July 2019.

\url{http://lazarus.elte.hu/cet/modules/guszlev/cylin.htm}}

\subsubsection*{Videos}

\item{Youtube, \textsl{Map projection of the Earth} (2011, December 12)

\url{https://www.youtube.com/watch?v=gGumy-9HrSY}}

\item{Vox, \textsl{Why all world maps are wrong?} (2016, December 2).

\url{https://www.youtube.com/watch?v=kIID5FDi2JQ}}

\item{Numberphile, \textsl{A strange map projection (Euler spiral)} (2018, November 13)

\url{https://www.youtube.com/watch?v=D3tdW9l1690}} 

\subsubsection*{Applications used}

\item{\LaTeX}

\item{Geogebra, \url{https://www.geogebra.org/}}

\item{Google maps}

\item{Derivative calculator, \url{https://www.derivative-calculator.net/}}

\item{Caliper - Mapping Software, GIS, and Transportation Software}

\item{Movable type scripts, \textsl{Calculate distance, bearing and more between latitude/longitude points} 

\url{https://www.movable-type.co.uk/scripts/latlong.html}}

\end{enumerate}

\end{document}